\documentclass[a4paper,fleqn]{cas-sc}

\usepackage[numbers]{natbib}

\usepackage{amssymb}
\usepackage{bm}

\usepackage{pgfplots}
\pgfplotsset{compat=1.18}
\usetikzlibrary{calc}
\usepackage{graphicx}
\usepackage{booktabs}    
\usepackage{multirow}    
\usepackage{graphicx}    
\usepackage{array}       
\usepackage{siunitx}     

\def\tsc#1{\csdef{#1}{\textsc{\lowercase{#1}}\xspace}}
\tsc{WGM}
\tsc{QE}

\begin{document}
\let\WriteBookmarks\relax
\def\floatpagepagefraction{1}
\def\textpagefraction{.001}

\shorttitle{Data-Boosted Optimization for AC-OPF}    

\shortauthors{Repiso et al.}

\title [mode = title]{Data-Boosted Optimization for AC Optimal Power Flow: Interior-Point and Spatial Branching Methods}  


\author[1]{Ignacio Repiso}[orcid=0009-0007-6386-776X]
\ead{irepiso@uma.es}
\credit{Writing – review \& editing, Writing – original draft, Visualization, Validation, Methodology, Investigation, Formal analysis, Data curation, Conceptualization}
\affiliation[1]{organization={OASYS Research Group, University of Malaga},
            city={Malaga},
            postcode={29071}, 
            state={Andalusia},
            country={Spain}}

\author[1]{Salvador Pineda}[orcid=0000-0002-1089-0970]
\ead{spineda@uma.es}
\credit{ Writing – review \& editing, Writing – original draft, Supervision, Project administration, Investigation, Funding acquisition, Conceptualization}

\author[1]{Juan Miguel Morales}[orcid=0000-0002-9114-686X]
\ead{juan.morales@uma.es}
\credit{Writing – review \& editing, Writing – original draft, Validation, Supervision, Investigation, Funding acquisition}

\begin{abstract}
The AC Optimal Power Flow (AC-OPF) problem is a non-convex, NP-hard optimization task essential for secure and economic power system operation. While interior-point methods are widely used due to their computational efficiency, spatial branching techniques offer global optimality guarantees at significantly higher computational cost.
In this work, we propose data-boosted variants of both approaches that leverage historical operating data to enhance performance. Specifically, data are used to guide initialization in interior-point methods and to restrict the search region in spatial branching. This unified perspective enables a systematic assessment of how learning can accelerate both local and global optimization strategies.
We conduct an extensive empirical study across networks of varying sizes under both standard conditions and modified configurations designed to induce local optima. Our results show that data-boosted strategies consistently improve convergence and reduce computation times for both approaches. However, spatial branching remains computationally demanding even with data-driven enhancements, while interior-point methods exhibit remarkable robustness, often converging to globally optimal solutions, even in challenging instances with multiple local optima.
These findings highlight the practical effectiveness of modern interior-point solvers and suggest that global optimization methods for AC-OPF still face significant scalability challenges, even when augmented with data-driven guidance.
\end{abstract}




\maketitle

\section{Introduction} \label{sec:intro}

The Alternating Current Optimal Power Flow (AC-OPF) problem is a fundamental tool for ensuring the secure, efficient, and economical operation of modern power systems. It determines the optimal operating point of the grid by minimizing generation costs subject to the nonlinear physical laws of power flow and numerous engineering constraints. However, the problem is inherently nonconvex, making globally optimal solutions computationally difficult to obtain \citep{sereeter2019optimal}. These challenges are compounded by the need to solve the problem repeatedly within short time intervals, as operators must rapidly adapt to demand fluctuations, renewable variability, and network contingencies \citep{tang2017real}. As an alternative to solving the full AC-OPF, the industry often relies on simplified or approximate models. A common example is the DC-OPF, which achieves computational efficiency through linearization, though at the expense of accuracy \citep{baker2021solutions}. Convex relaxations offer another widely studied approach, providing theoretical guarantees of global optimality; however, they can yield solutions that are infeasible for the original nonconvex problem and frequently entail considerable computational effort \citep{molzahn2019survey}.

To guarantee physically feasible solutions, the AC-OPF has traditionally been addressed using interior-point methods and other nonlinear optimization techniques. These approaches are generally efficient in practice but may converge to local rather than global optima \citep{narimani2018empirical} or even face convergence difficulties \citep{capitanescu2013experiments}. One way to mitigate this drawback is through multiple solver initializations, although this substantially increases computational effort \citep{zhang2024learning}. Despite the non-convex nature of the AC-OPF, recent studies show that interior-point algorithms often converge to the global optimum in practice. For example, \citep{kardos2018complete} reports that, across a wide range of networks and operating conditions, different initializations typically lead to the same globally optimal solution. Similarly, the exhaustive experiments in \citep{bukhsh2013local} indicate that even for cases deliberately designed to exhibit local optima, interior-point methods still recover the global solution, including when initialized with the flat voltage profile.

More recently, spatial branching has emerged as an alternative methodology to classical interior-point solvers for the AC optimal power flow problem. In this approach, the problem is reformulated as a quadratically constrained quadratic program (QCQP) in rectangular voltage coordinates and solved through a systematic partitioning of the variable space combined with bound-tightening techniques \citep{chen2015bound}. Further developments have introduced more sophisticated branching strategies to prioritize the most promising variables \citep{iranpour2025scalable}. In principle, these methods can guarantee global optimality and feasibility, which makes them an appealing theoretical alternative to interior-point algorithms. However, their computational requirements remain extremely demanding, even for networks of moderate size, which raises questions about their practical applicability.


Both interior-point and spatial branching methods rely heavily on problem-specific parameters that critically influence their performance. For interior-point algorithms, the quality of the obtained solution depends on the initial operating point, while for spatial branching, computational efficiency is largely determined by the variable bounds defining the search region. Motivated by recent trends in learning-based optimization \citep{lei2020data, stratigakos2023interpretable, pineda2025beyond}, we introduce data-boosted variants of both strategies. In the interior-point approach, historical system data are used to identify the most similar past operating condition and initialize the solver at its optimal solution. In spatial branching, a $K$-nearest neighbors ($K$NN) strategy learns adjusted voltage bounds from historical data, effectively reducing the search space. These two data-boosted variants leverage past operational knowledge to improve convergence and computational performance while preserving feasibility and interpretability. According to \citep{khaloie2025review}, machine learning methods for AC-OPF are commonly classified into \emph{End-to-End} (E2E) approaches, which predict solutions directly from system inputs but may violate constraints and therefore require post-processing to restore feasibility \citep{pan2020deepopf,taheri2024ac}, and \emph{Learning-to-Optimize} (L2O) methods, which integrate learning with classical optimization to accelerate or guide the solution process \citep{cengil2022learning,crozier2022data}. Our data-boosted approaches clearly fall within the L2O paradigm, combining learning with optimization to enhance efficiency.

The primary contribution of this paper is the development of a data-boosted optimization framework for AC-OPF that integrates learning into both interior-point and spatial branching methods. By leveraging historical operating data to guide initialization and restrict the search space, the proposed approach provides a unified perspective on how data can enhance both local and global optimization strategies. Building on this framework, we conduct an extensive empirical study across a diverse set of networks and operating conditions, including both standard benchmark cases and challenging configurations designed to induce local optima. This analysis reveals fundamental differences in the practical behavior of these methods and provides new insights into their performance, robustness, and scalability under varying levels of problem difficulty.


The remainder of the paper is organized as follows. Section~\ref{sec:formulation} introduces the nomenclature and presents the formulation of the AC optimal power flow problem. Section~\ref{sec:methodology} details the data-boosted versions of both the interior-point and spatial-branching approaches. The performance evaluation metrics used to compare these methods are described in Section~\ref{sec:comparison}. Section~\ref{sec:computational_results} reports and discusses the computational results obtained for networks of different sizes and operating conditions. Finally, Section~\ref{sec:conclusions} summarizes the main findings and outlines the conclusions of the study.

\section{Formulation} \label{sec:formulation}

We consider a power system consisting of nodes, branches interconnecting them, and generating units. Each node, indexed by $n \in \mathcal{N}$, is described by its shunt conductance ($G^{sh}_n$), shunt susceptance ($B^{sh}_n$), active load ($p^d_n$), and reactive load ($q^d_n$). The state variables at each node are now represented in rectangular coordinates as $e_n$ and $f_n$, corresponding to the real and imaginary components of the voltage $V_n = e_n + \mathrm{j}f_n$. The voltage magnitude is implicitly constrained by $\underline{v}_n^2 \le e_n^2 + f_n^2 \le \overline{v}_n^2$. Generators, indexed by $g \in \mathcal{G}$, are associated with convex generation cost functions $c_g(\cdot)$ and must respect technical limits on their outputs, i.e., $\underline{p}_g \leq p_g \leq \overline{p}_g$ for active power and $\underline{q}_g \leq q_g \leq \overline{q}_g$ for reactive power. Branches, indexed by $l \in \mathcal{L}$, are characterized by electrical parameters including resistance $r_l$, reactance $x_l$, total charging susceptance $b_l$, tap ratio $\tau_l$, and phase shift $\gamma_l$, all assumed to be known. To capture network connectivity, we define two sparse matrices $F$ and $T$ of dimension $|\mathcal{L}| \times |\mathcal{N}|$, where $F_{ln}$ (resp. $T_{lm}$) equals 1 if node $n$ ($m$) is the sending (receiving) end of branch $l$, and 0 otherwise. For each branch $l$, the active power sent and received are denoted $p^f_l$ and $p^t_l$, respectively, while $q^f_l$ and $q^t_l$ denote the corresponding reactive flows. Branch apparent power is constrained by the thermal capacity limit $\overline{s}_l$. Finally, the branch-end admittance coefficients are given by
\begin{subequations}
\begin{align}
& G^{ff}_l + jB^{ff}_l = \left(\tfrac{1}{r_l+jx_l}+j\tfrac{b_l}{2} \right)\tfrac{1}{\tau_l^2}, \nonumber\\
& G^{ft}_l + jB^{ft}_l = -\tfrac{1}{r_l+jx_l}\tfrac{1}{\tau_l e^{-j\gamma_l}}, \nonumber\\
& G^{tf}_l + jB^{tf}_l = -\tfrac{1}{r_l+jx_l}\tfrac{1}{\tau_l e^{j\gamma_l}}, \nonumber\\
& G^{tt}_l + jB^{tt}_l = \left(\tfrac{1}{r_l+jx_l}+j\tfrac{b_l}{2} \right) \nonumber.
\end{align}
\end{subequations}

Under this notation, the AC-OPF problem can be formulated as the following optimization model \citep{nair2022computational, bienstock2020mathematical}:
\begin{subequations}\label{eq:ac_opf_rect}
\begin{align}
& \min_{p_g,p^f_l,q^f_l,p^t_l,q^t_l,e_n,f_n} \sum_g c_g(p_g),  \label{eq:ac_opf_rect_of}\\
& \text{s.t.} \nonumber\\
& \sum_{g\in\mathcal{G}_n} p_g - p^d_n = (e_n^2+f_n^2)G^{sh}_n + \sum_{l} (F_{ln}p^f_l + T_{ln}p^t_l), &\forall n,\label{eq:ac_opf_rect_balp}\\
& \sum_{g\in\mathcal{G}_n} q_g - q^d_n = -(e_n^2+f_n^2)B^{sh}_n + \sum_{l} (F_{ln}q^f_l + T_{ln}q^t_l), &\forall n,\label{eq:ac_opf_rect_balq}\\
\displaybreak[3] 
& p^f_l = G^{ff}_l (e_n^2+f_n^2) + G^{ft}_l (e_n e_m + f_n f_m)  \nonumber\\
&\qquad - B^{ft}_l (e_n f_m - f_n e_m), \quad \forall(l,n,m):F_{ln}=T_{lm}=1,\label{eq:ac_opf_rect_pf}\\
& q^f_l = -B^{ff}_l (e_n^2+f_n^2) - B^{ft}_l (e_n e_m + f_n f_m)\nonumber\\
&\qquad - G^{ft}_l (e_n f_m - f_n e_m), \quad \forall(l,n,m):F_{ln}=T_{lm}=1,\label{eq:ac_opf_rect_qf}\\
\displaybreak[3]
& p^t_l = G^{tt}_l (e_m^2+f_m^2) + G^{tf}_l (e_m e_n + f_m f_n) \nonumber\\
&\qquad - B^{tf}_l (e_m f_n - f_m e_n), \quad \forall(l,n,m):F_{ln}=T_{lm}=1,\label{eq:ac_opf_rect_pt}\\
& q^t_l = -B^{tt}_l (e_m^2+f_m^2) - B^{tf}_l (e_m e_n + f_m f_n)\nonumber\\
&\qquad - G^{tf}_l (e_m f_n - f_m e_n), \quad \forall(l,n,m):F_{ln}=T_{lm}=1,\label{eq:ac_opf_rect_qt}\\
\displaybreak[3]
& \underline{p}_g \leq p_g \leq \overline{p}_g, &\forall g, \label{eq:ac_opf_rect_pglim}\\
& \underline{q}_g \leq q_g \leq \overline{q}_g, &\forall g, \label{eq:ac_opf_rect_qglim}\\
& \underline{v}_n^2 \leq e_n^2 + f_n^2 \leq \overline{v}_n^2, &\forall n, \label{eq:ac_opf_rect_vlim}\\
& (p^f_l)^2 + (q^f_l)^2 \leq (\overline{s}_l)^2, &\forall l, \label{eq:ac_opf_rect_sfmax}\\
& (p^t_l)^2 + (q^t_l)^2 \leq (\overline{s}_l)^2, &\forall l, \label{eq:ac_opf_rect_stmax}
\end{align}
\end{subequations}

\noindent \noindent where $\mathcal{G}_n$ denotes the set of generators connected to node $n$. The objective function~\eqref{eq:ac_opf_rect_of} aims to minimize the total generation cost of the system. Node-level active and reactive power balances are enforced by~\eqref{eq:ac_opf_rect_balp} and~\eqref{eq:ac_opf_rect_balq}, respectively. The active and reactive flows on each branch $l$ are described by the equalities~\eqref{eq:ac_opf_rect_pf}--\eqref{eq:ac_opf_rect_qt}. Operational and technical limits are imposed by constraints~\eqref{eq:ac_opf_rect_pglim}--\eqref{eq:ac_opf_rect_stmax}. Although the cost function~\eqref{eq:ac_opf_rect_of} is typically quadratic, linear, or piecewise linear, the overall formulation~\eqref{eq:ac_opf_rect} remains non-convex due to the quadratic and bilinear terms in the rectangular voltage variables appearing in~\eqref{eq:ac_opf_rect_pf}--\eqref{eq:ac_opf_rect_qt}.

\section{Methodology} \label{sec:methodology}

In this section, we describe the approaches used to solve the non-convex optimization model~\eqref{eq:ac_opf_rect}, focusing on two state-of-the-art strategies: interior-point algorithms and spatial branching. For each strategy, we consider standard versions that solve the problem using only the current instance’s data, as well as data-boosted variants that leverage historical solutions to enhance performance. In the interior-point case, historical data are used to initialize the solver at the optimal solution of the most similar past operating point, while for spatial branching, a $K$-nearest neighbors approach adjusts the variable bounds to reduce the search space. For these data-boosted variants, we assume access to a dataset of historical solutions obtained offline by solving model~\eqref{eq:ac_opf_rect} to global optimality. Each entry in this dataset consists of an operating point characterized by the profiles of active and reactive power demand, $\bm{p}^d$ and $\bm{q}^d$, together with the corresponding optimal solution of the AC-OPF, represented by the real and imaginary components of the bus voltages, $\bm{e}^*$ and $\bm{f}^*$, respectively. Formally, this can be expressed as:
\[
\{(\bm{p}^d_i, \bm{q}^d_i, \bm{e}^*_i, \bm{f}^*_i) \,|\, i \in \mathcal{I}\},
\]
where $\mathcal{I}$ denotes the index set of training instances.

We first present the methodology that solves the AC-OPF using interior-point (IP) methods, which have historically been the most widely adopted approach. Starting from an initial voltage profile, the method iteratively refines the solution by applying Newton-type steps to a barrier reformulation of the problem \citep{ye2011interior}. This approach is computationally efficient and scales well to large networks, but due to the non-convexity of the AC-OPF, it may converge only to local optima. In our study, we investigate several variants of the interior-point approach. Prior works \citep{kardos2018complete} have shown that OPF performance depends not only on the chosen solution method, but also on the mathematical formulation used. To account for this, we solve the AC-OPF in both rectangular (formulated in \eqref{eq:ac_opf_rect}) and polar coordinates (formulated in reference \citep{sereeter2019optimal}). We consider two initialization strategies for the interior-point solver. The first is the flat voltage profile, which is the most commonly used in practice. In this case, the initial voltage components are set uniformly across all nodes as
\[
e^0_n = 1, \quad f^0_n = 0, \quad \forall n \in \mathcal{N}.
\]
The second strategy is a historical initialization, corresponding to our data-boosted variant. Here, the solver is initialized using the optimal solution of the historical instance whose demand profile is closest to the current one. Formally, the initial point is defined as
\[
\bm{e}^0 = \bm{e}^*_{i'}, \quad 
\bm{f}^0 = \bm{f}^*_{i'}, \quad 
i' \in \arg\min_{i \in \mathcal{I}} \, ||\bm{p}^d - \bm{p}^d_i||^2_2 +  ||\bm{q}^d - \bm{q}^d_i||^2_2 ,
\]
where \(\mathcal{I}\) denotes the set of historical instances.


We next consider the use of spatial branching (SP) techniques to solve the AC-OPF formulated as a quadratically constrained quadratic program (QCQP) in rectangular voltage coordinates. Spatial branching systematically partitions the feasible region into smaller subregions and solves convex relaxations in each, progressively refining bounds until convergence to the global optimum is achieved \citep{linderoth2005simplicial}. While this approach guarantees global optimality, its computational burden can still be significant in practice.

\begin{figure}
    \centering
    \includegraphics[width=0.45\linewidth]{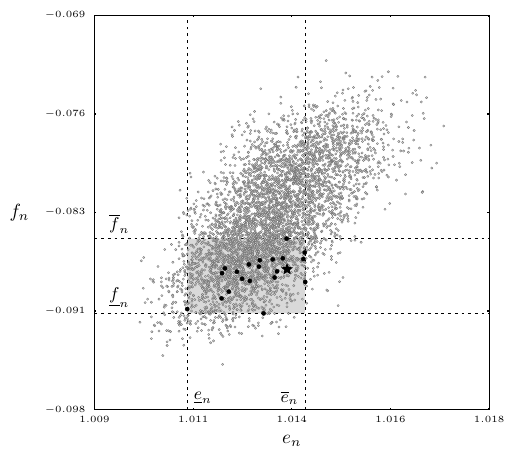}
    \caption{Voltage bounds reduction for \texttt{SP-\scriptsize{$K$}} approach.}
    \label{fig:learned_bounds}
\end{figure}

To assess its performance, we consider three variations. The first, SP--O, uses the original bounds of the voltage variables computed as follows. Given the voltage magnitude limits $\underline v_n$ and $\overline v_n$ and the voltage angle bounds $\underline\theta_n$ and $\overline\theta_n$, we first define candidate sets for the voltage magnitude and angle:
\[
\begin{aligned}
\mathcal V_n &= \{\underline v_n,\, \overline v_n\}, \quad &
    \forall n &\in \mathcal{N}, \\
\mathcal A_n &= \{\theta \in \{-\pi,-\tfrac{\pi}{2},0,\tfrac{\pi}{2},\pi\}
                \mid \underline\theta_n \le \theta \le \overline\theta_n \}
                \,\cup\, \{\underline\theta_n,\,\overline\theta_n\} \quad &
    \forall n &\in \mathcal{N},
\end{aligned}
\]
\noindent where, $\mathcal V_n$ contains the extreme values of the voltage magnitude, while $\mathcal A_n$ includes selected canonical angles within the allowed range as well as the angle bounds themselves. These candidate sets are then used to determine conservative bounds on the rectangular voltage components by evaluating all combinations of voltage magnitude and angle:
\begin{subequations}
\label{eq:bounds_rect}
\begin{alignat*}{4}
    \underline e^0_n &= \min_{v\in\mathcal V_n,\,\theta\in\mathcal A_n} v\cos\theta, \quad &
    \overline e^0_n  &= \max_{v\in\mathcal V_n,\,\theta\in\mathcal A_n} v\cos\theta, \quad &
    \forall n &\in \mathcal{N}, \\
    \underline f^0_n &= \min_{v\in\mathcal V_n,\,\theta\in\mathcal A_n} v\sin\theta, \quad &
    \overline f^0_n  &= \max_{v\in\mathcal V_n,\,\theta\in\mathcal A_n} v\sin\theta, \quad &
    \forall n &\in \mathcal{N}.
\end{alignat*}
\end{subequations}
These bounds provide a rectangular envelope that fully contains all feasible voltage values in the polar representation, serving as the starting point for the spatial branching procedure. Second, we consider SP-$\epsilon$, a benchmark variant in which the bounds of the voltage variables are artificially tightened around the known optimal solution:
\begin{alignat*}{4}
    \underline{e}_n &= \max(\underline{e}^0_n, e^*_n - \epsilon), \quad &
    \overline{e}_n  &= \min(\overline{e}^0_n, e^*_n + \epsilon), \quad &
    \forall n &\in \mathcal{N}, \\
    \underline{f}_n &= \max(\underline{f}^0_n, f^*_n - \epsilon), \quad &
    \overline{f}_n  &= \min(\overline{f}^0_n, f^*_n + \epsilon), \quad &
    \forall n &\in \mathcal{N}.
\end{alignat*}
where $\epsilon > 0$ is a small positive constant. This idealized construction assumes perfect knowledge of the optimal solution, and therefore provides a lower bound on the computational effort that could be achieved by any bounding or search-space reduction strategy. Third, we consider a data-boosted variant, denoted SP--{\footnotesize{$K$}}, where $K$ specifies the number of nearest neighbors selected from a historical dataset. In this approach, we first identify the $K$ historical operating points whose demand profiles are closest to the current one, based on the Euclidean distance $||\bm{p}^d - \bm{p}^d_i||^2_2 +  ||\bm{q}^d - \bm{q}^d_i||^2_2$. We denote the set of corresponding indices as $\mathcal{I}^K$. The bounds on the voltage variables are then computed from these neighbors as
\begin{alignat*}{4}
    \underline{e}_n &= \min_{i\in\mathcal{I}^K} e^*_{in}, \quad &
    \overline{e}_n  &= \max_{i\in\mathcal{I}^K} e^*_{in}, \quad &
    \forall n &\in \mathcal{N}, \\
    \underline{f}_n &= \min_{i\in\mathcal{I}^K} f^*_{in}, \quad &
    \overline{f}_n  &= \max_{i\in\mathcal{I}^K} f^*_{in}, \quad &
    \forall n &\in \mathcal{N}.
\end{alignat*}
\noindent where $e^*_{in}$ and $f^*_{in}$ denote the real and imaginary components, respectively, of the optimal voltage at bus $n$ for instance $i$. In addition, the data-driven bounds are rounded up or down to the fifth decimal place to slightly enlarge the feasible region and prevent potential numerical issues related to the solver's feasibility tolerances. Figure~\ref{fig:learned_bounds} illustrates, in the complex plane, the learned voltage bounds obtained using the data-boosted \texttt{SP-\footnotesize{$K$}} approach for a specific bus. Gray dots represent the optimal complex voltages of this bus across all training instances. The star marks the optimal solution of a particular test instance, while the bold dots indicate the $K$ nearest neighbors of this test instance selected from the historical dataset. The dashed lines highlight the tightened voltage bounds derived from these neighbors, which are the key element used to restrict the search space for the spatial branching algorithm. By constraining the feasible region to the bounds indicated by the dashed lines, the data-boosted approach can help alleviate the computational burden of spatial branching. The main limitation of this heuristic strategy is that it cannot theoretically guarantee that the global optimum of the original problem lies within the restricted bounds. Nevertheless, as shown in the figure, the optimal solution for this test instance falls within the tightened region. Across our experiments, the probability of excluding the true global optimum was negligible and occurred only in the most challenging instances prone to local optima.

\section{Comparison} \label{sec:comparison}

In this section, we describe the setup and methodology used to compare the different approaches for solving the AC-OPF problem. Our analysis encompasses both interior-point and spatial branching methods, considering their standard versions as well as their data-boosted variants. For clarity, the specific configurations evaluated are summarized below:

\begin{itemize}
    \item \texttt{IP-R-F}: Interior-point solver using the \emph{rectangular} formulation, initialized with a flat voltage profile ($e_n=1$, $f_n=0$ for all buses).
    \item \texttt{IP-R-H}: Interior-point solver using the rectangular formulation, initialized with the nearest historical operating point (data-boosted variant).
    \item \texttt{IP-P-F}: Interior-point solver using the \emph{polar} formulation, initialized with a \emph{flat voltage profile}, i.e., $|V_n|=1$ p.u. and $\theta_n=0$ for all buses.
    \item \texttt{IP-P-H}: Interior-point solver using the polar formulation, initialized with the optimal solution of the most similar historical demand instance (data-boosted variant).
    \item \texttt{SP-O}: Spatial branching solver using the original voltage bounds derived from the operational limits of voltage magnitude and phase angles.
    \item \texttt{SP-$\epsilon$}: Benchmark spatial branching variant with artificially tightened voltage bounds around the optimal solution (idealized lower-bound case).
    \item \texttt{SP-\footnotesize{$K$}}: Data-boosted spatial branching variant in which voltage bounds are learned from the $K$ nearest historical operating points, effectively reducing the search region.
\end{itemize}

We evaluate the performance of all methods using four complementary metrics. First, the \emph{computational time} ($T$), measured in seconds, quantifies the efficiency of each approach. Second, the \emph{suboptimality} (\texttt{sub}) measures the relative deviation of the obtained objective value from the global optimum (when available) or, otherwise, from the best solution among all approaches. Third, for the spatial-branching methods, we report the \emph{MIP gap} (\texttt{gap}), defined as the relative difference between the best incumbent solution and the best bound at termination. Finally, we include the number of infeasible instances, denoted by \texttt{\#inf}, which is a particularly relevant metric for assessing the performance of spatial-branching methods, since in some cases the solver fails to find a feasible solution within the time limit. Taken together, these four metrics provide a comprehensive evaluation of computational efficiency, solution quality, and optimality.

\section{Computational Results} \label{sec:computational_results}

In this section, we evaluate the proposed data-boosted framework described in Section~\ref{sec:methodology}, using four IEEE test cases with 14, 39, 57, and 118 buses. For each network, we analyze two different operating-condition settings. The first corresponds to the original data available in Matpower~\cite{zimmerman2010matpower}, where non-convexities are mild and the global optimum can be obtained reliably. The second setting introduces modifications to the network data, making the presence of local optima more likely and the problem more challenging. Specifically, for the 39- and 118-bus networks we adopt the modifications proposed in~\cite{bukhsh2013local}, while for the 14- and 57-bus networks we apply the modifications proposed in~\cite{narimani2018empirical}. Detailed network and generator specifications, together with the nominal demand profiles, can be found in~\cite{Oasys2025repository}.

For each of the eight network types and configurations, we generate a dataset used by the learning-based approaches that rely on historical data. The construction procedure is as follows. Let $\hat{d}_n$ denote the nominal demand at node $n$ as provided in the repository~\cite{Oasys2025repository}. We create 5\,000 demand profiles by sampling each nodal demand from a uniform distribution within $\pm 5\%$ of its nominal value, i.e.,
\[
d_n \sim \hat{d}_n \cdot \mathcal{U}([0.95,1.05]).
\] 
For each instance $i$, defined by a sampled demand profile, we solve the AC-OPF problem using three baseline approaches that do not rely on historical data: \texttt{IP-P-F}, \texttt{IP-R-F}, and \texttt{SP-O}. The two interior-point formulations (\texttt{IP-P-F} and \texttt{IP-R-F}) always provide a feasible solution within the one-hour time limit, although it may correspond to a local optimum. In contrast, the spatial branching approach (\texttt{SP-O}) may require more time, and in some cases, it does not converge within one hour. However, whenever \texttt{SP-O} finishes, the obtained solution is guaranteed to be globally optimal because the solver certifies optimality through the MIP gap. For each instance, we first check the feasibility of the solutions provided by the three methods, retaining only those with a maximum constraint violation below $10^{-5}$. Among the feasible ones, we select the solution with the lowest objective value. If \texttt{SP-O} provides a certified global optimum, that solution is used. Otherwise, the best feasible solution among the interior-point methods is selected. The resulting solution is then stored in the dataset of historical OPF results.

From the 5\,000 instances, we randomly assign 4\,900 to the training set and 100 to the test set. The training set is used to learn the initialization points and variable bounds for the data-boosted approaches, namely \texttt{IP-P-H}, \texttt{IP-R-H}, \texttt{SP-$\epsilon$}, and \texttt{SP-\footnotesize{$K$}}. The test set is reserved for evaluation: all approaches (both baseline and data-boosted) are independently applied to these 100 unseen instances, and the reported results correspond to their performance on this common benchmark. For the \texttt{SP-$\epsilon$} approach, which restricts the search region to a small neighborhood around the optimal solution, the value of $\epsilon$ is set to $10^{-5}$. Finally, to analyze the sensitivity of the \texttt{SP-\footnotesize{$K$}} method, we consider two neighborhood sizes, $K \in \{$20, 100\}. Smaller values of $K$ lead to tighter learned bounds and thus have greater potential to accelerate the spatial branching process, although at the risk of excluding the global optimum from the reduced feasible region. 

All optimization problems are formulated in \texttt{AMPL} and solved using \texttt{IPOPT~3.12.13} for the interior-point approaches and \texttt{Gurobi~12.0.2} for the spatial branching approaches. Simulations are executed on a Linux-based server with an AMD~EPYC processor running at 2.25~GHz, using a single thread and 8~GB of~RAM. In all experiments, the optimality gap tolerance is set to $0.01\%$, and the maximum computation time is limited to one hour.

In Table~\ref{tab:double_results}, we summarize the computational results obtained for the eight network configurations and all considered approaches. The reported values correspond to averages over the 100 test instances for computational time~($\widehat{T}$), suboptimality~($\widehat{\texttt{sub}}$), and MIP gap~($\widehat{\texttt{gap}}$). In addition, we also report the maximum values observed among the 100 instances, indicated with a hat as $\overline{T}$, $\overline{\texttt{sub}}$, and $\overline{\texttt{gap}}$. If the computational time is below one second, it is reported as $<1$, whereas if the 1-hour time limit is reached, it is indicated as TL. For the spatial-branching approaches, a check mark ($\checkmark$) denotes that the final MIP gap for all instances is below the threshold of 0.01\%. For the interior-point approaches, the MIP gap is not applicable and is represented by a dash (-) in the table. We also include the number of infeasible cases, indicated as $\texttt{\#inf}$.



\begin{table*}[h]
\caption{Summary of computational results across network sizes, operating conditions, and AC-OPF solution approaches.}
\renewcommand{\arraystretch}{1.2}
\centering
\resizebox{\textwidth}{!}{
\begin{tabular}{l l cc cc cc c | cc cc cc c}
\toprule
& & \multicolumn{7}{c|}{\footnotesize \textbf{Standard configurations}} &
    \multicolumn{7}{c}{\footnotesize \textbf{Tuned configurations}} \\
\cmidrule(lr){3-9} \cmidrule(lr){10-16}
\textbf{Test case} & \textbf{Approach} 
  & $\widehat{T}$ (s) & $\overline{T}$ (s) 
  & $\widehat{\texttt{gap}}$ (\%) & $\overline{\texttt{gap}}$ (\%) 
  & $\widehat{\texttt{sub}}$ (\%) & $\overline{\texttt{sub}}$ (\%) 
  & \texttt{\#inf} 
  & $\widehat{T}$ (s) & $\overline{T}$ (s) 
  & $\widehat{\texttt{gap}}$ (\%) & $\overline{\texttt{gap}}$ (\%) 
  & $\widehat{\texttt{sub}}$ (\%) & $\overline{\texttt{sub}}$ (\%) 
  & \texttt{\#inf} \\
\midrule
\multirow{9}{*}{IEEE 14} & \texttt{IP-R-F} & $<$1 & $<$1 & - & - & 1.0e-3 & 2.0e-3 & 0 & $<$1 & $<$1 & - & - & 2.2e-3 & 4.7e-3 & 0 \\
 & \texttt{IP-R-H} & $<$1 & $<$1 & - & - & 1.0e-3 & 2.0e-3 & 0 & $<$1 & $<$1 & - & - & 2.2e-3 & 4.7e-3 & 0 \\
 & \texttt{IP-P-F} & $<$1 & $<$1 & - & - & 1.0e-3 & 2.0e-3 & 0 & $<$1 & $<$1 & - & - & 2.2e-3 & 4.7e-3 & 0 \\
 & \texttt{IP-P-H} & $<$1 & $<$1 & - & - & 1.0e-3 & 2.0e-3 & 0 & $<$1 & $<$1 & - & - & 2.2e-3 & 4.7e-3 & 0 \\
 & \texttt{SP-O} & 24 & 1149 & $\checkmark$ & $\checkmark$ & 1.0e-3 & 3.0e-3 & 0 & 29 & 192 & $\checkmark$ & $\checkmark$ & 9.1e-4 & 3.9e-3 & 0 \\
 & \texttt{SP-$\epsilon$} & $<$1 & $<$1 & $\checkmark$ & $\checkmark$ & 1.0e-3 & 2.2e-3 & 0 & $<$1 & $<$1 & $\checkmark$ & $\checkmark$ & 2.4e-3 & 5.2e-3 & 0 \\
 & \texttt{SP-20} & $<$1 & $<$1 & $\checkmark$ & $\checkmark$ & 1.0e-3 & 3.5e-3 & 1 & $<$1 & 1.2 & $\checkmark$ & $\checkmark$ & 1.5e-3 & 8.2e-3 & 3 \\
 & \texttt{SP-100} & $<$1 & $<$1 & $\checkmark$ & $\checkmark$ & 3.7e-4 & 1.8e-3 & 0 & $<$1 & $<$1 & $\checkmark$ & $\checkmark$ & 2.2e-3 & 4.7e-3 & 0 \\
\midrule
\multirow{9}{*}{IEEE 39} & \texttt{IP-R-F} & $<$1 & $<$1 & - & - & 3.8e-4 & 8.6e-4 & 0 & $<$1 & $<$1 & - & - & 2.9e-2 & 5.0e-2 & 0 \\
 & \texttt{IP-R-H} & $<$1 & 1.3 & - & - & 3.8e-4 & 8.6e-4 & 0 & $<$1 & $<$1 & - & - & 2.9e-2 & 5.0e-2 & 0 \\
 & \texttt{IP-P-F} & $<$1 & $<$1 & - & - & 3.8e-4 & 8.6e-4 & 0 & $<$1 & $<$1 & - & - & 2.9e-2 & 5.0e-2 & 0 \\
 & \texttt{IP-P-H} & $<$1 & $<$1 & - & - & 3.8e-4 & 8.6e-4 & 0 & $<$1 & $<$1 & - & - & 2.9e-2 & 5.0e-2 & 0 \\
 & \texttt{SP-O} & 2910 & TL & 1.1e-2 & 2.2e-2 & 3.7e-4 & 9.2e-4 & 0 & TL & TL & 3.7e+0 & 7.4e+0 & 1.3e-2 & 4.1e-2 & 0 \\
 & \texttt{SP-$\epsilon$} & $<$1 & $<$1 & $\checkmark$ & $\checkmark$ & 4.3e-4 & 8.6e-4 & 0 & $<$1 & $<$1 & $\checkmark$ & $\checkmark$ & 1.1e-2 & 3.9e-2 & 0 \\
 & \texttt{SP-20} & 1.7 & 3.9 & $\checkmark$ & $\checkmark$ & 4.0e-4 & 5.2e-3 & 0 & 3562 & TL & 6.4e-1 & 1.0e+1 & 1.9e-1 & 1.0e+1 & 2 \\
 & \texttt{SP-100} & 1.5 & 2.4 & $\checkmark$ & $\checkmark$ & 2.5e-4 & 7.0e-4 & 0 & TL & TL & 9.6e-1 & 1.9e+0 & 1.7e-2 & 4.1e-2 & 2 \\
\midrule
\multirow{9}{*}{IEEE 57} & \texttt{IP-R-F} & $<$1 & $<$1 & - & - & 1.4e-3 & 3.3e-3 & 0 & 1.4 & 6.1 & - & - & 5.4e-1 & 1.5e+0 & 2 \\
 & \texttt{IP-R-H} & $<$1 & 1.0 & - & - & 1.4e-3 & 3.3e-3 & 0 & $<$1 & $<$1 & - & - & 1.7e-1 & 1.3e+0 & 0 \\
 & \texttt{IP-P-F} & $<$1 & $<$1 & - & - & 1.4e-3 & 3.3e-3 & 0 & $<$1 & 2.2 & - & - & 5.1e-1 & 1.5e+0 & 0 \\
 & \texttt{IP-P-H} & $<$1 & $<$1 & - & - & 1.4e-3 & 3.3e-3 & 0 & $<$1 & $<$1 & - & - & 1.6e-1 & 1.3e+0 & 0 \\
 & \texttt{SP-O} & TL & TL & 6.2e-2 & 6.6e-2 & 1.4e-3 & 3.5e-3 & 0 & TL & TL & 2.3e+1 & 3.8e+1 & 3.9e-2 & 1.8e-1 & 59 \\
 & \texttt{SP-$\epsilon$} & $<$1 & $<$1 & $\checkmark$ & $\checkmark$ & 6.8e-4 & 2.7e-3 & 0 & $<$1 & 2.0 & $\checkmark$ & $\checkmark$ & 2.0e-1 & 1.1e+0 & 1 \\
 & \texttt{SP-20} & 5.5 & 56 & $\checkmark$ & $\checkmark$ & 1.2e-3 & 4.6e-3 & 23 & 3562 & TL & 1.9e+1 & 4.7e+1 & 9.0e-2 & 1.7e+0 & 50 \\
 & \texttt{SP-100} & 6.0 & 73 & $\checkmark$ & $\checkmark$ & 1.3e-3 & 3.3e-3 & 0 & 3577 & TL & 1.6e+1 & 4.6e+1 & 5.9e-2 & 7.5e-1 & 39 \\
\midrule
\multirow{9}{*}{IEEE 118} & \texttt{IP-R-F} & 2.7 & 8.7 & - & - & 4.0e-4 & 1.1e-3 & 0 & 3.4 & 6.0 & - & - & 4.9e-4 & 1.1e-3 & 0 \\
 & \texttt{IP-R-H} & 1.8 & 3.2 & - & - & 4.0e-4 & 1.1e-3 & 0 & 2.6 & 5.8 & - & - & 4.9e-4 & 1.1e-3 & 0 \\
 & \texttt{IP-P-F} & $<$1 & $<$1 & - & - & 4.0e-4 & 1.1e-3 & 0 & $<$1 & $<$1 & - & - & 4.9e-4 & 1.1e-3 & 0 \\
 & \texttt{IP-P-H} & $<$1 & $<$1 & - & - & 4.0e-4 & 1.1e-3 & 0 & $<$1 & $<$1 & - & - & 4.9e-4 & 1.1e-3 & 0 \\
 & \texttt{SP-O} & TL & TL & 2.8e-1 & 3.0e-1 & 3.7e-4 & 1.2e-3 & 0 & TL & TL & 2.7e-1 & 2.9e-1 & 4.6e-4 & 1.2e-3 & 0 \\
 & \texttt{SP-$\epsilon$} & 1.1 & 1.7 & $\checkmark$ & $\checkmark$ & 4.0e-4 & 1.1e-3 & 0 & 1.0 & 1.7 & $\checkmark$ & $\checkmark$ & 5.0e-4 & 1.1e-3 & 0 \\
 & \texttt{SP-20} & 15 & 398 & $\checkmark$ & $\checkmark$ & 5.9e-4 & 1.4e-2 & 1 & 11 & 21 & $\checkmark$ & $\checkmark$ & 2.9e-4 & 1.1e-3 & 1 \\
 & \texttt{SP-100} & 16 & 27 & $\checkmark$ & $\checkmark$ & 3.5e-4 & 1.1e-3 & 0 & 14 & 26 & $\checkmark$ & $\checkmark$ & 4.3e-4 & 1.1e-3 & 0 \\
\bottomrule
\end{tabular}
}
\label{tab:double_results}
\end{table*}

We begin by analyzing the results for the two 14-bus networks. Owing to their small size, all interior-point approaches solve the OPF in less than one second and achieve solutions that are effectively globally optimal, as reflected by the very low \texttt{sub} values. In contrast, even for this modest network size, the \texttt{SP-O} approach requires substantially more computational time, although it provides a certificate of global optimality. The data-boosted spatial-branching variants markedly reduce computational time, achieving speedups of up to 288$\times$ over \texttt{SP-O}, without compromising optimality as indicated by the low $\texttt{sub}$ values. However, for some test instances, the tightened bounds lead to infeasibility. This issue is mitigated by increasing the number of neighbors from $K=$ 20 to $K=$ 100, which slightly enlarges the learned bounds and ensures feasibility across all 100 test instances while still maintaining considerable reductions in computational effort. No substantial differences are observed between the standard and tuned versions of the 14-bus network.

For the two 39-bus networks, the solution times of the interior-point methods increase slightly but remain low. Although this system size is still moderate, the spatial-branching approach with default bounds \texttt{SP-O} reaches the one-hour time limit in some instances for the standard configuration and in all instances for the tuned configuration, which is more prone to local optima. In both cases, the average $\texttt{gap}$ exceeds 0.01\%, indicating that global optimality could not be certified within the allotted time. Nevertheless, the small $\texttt{sub}$ values suggest that the incumbent solution is very close to, and likely coincides with, the true global optimum in most cases, illustrating a key limitation of the spatial-branching approach, which often fails to certify optimality despite finding near-optimal solutions. For the standard 39-bus network, the data-boosted spatial-branching methods substantially reduce computational time, allowing all instances to be solved in under 10 seconds while still attaining globally optimal solutions. For the tuned version, which exhibits a more complex and multimodal landscape, the data-boosted variants reduce the average $\texttt{gap}$ but still hit the time limit in most cases. Some instances also become infeasible due to overly tightened bounds, meaning that the true optimal solution lies outside the restricted search region. This issue is evident in one test case where the maximum $\texttt{sub}$ reaches 10\%, a problem that disappears when increasing $K$ from 20 to 100.

For the 57-bus networks, the interior-point approaches continue to show reasonable computational times. However, the tuned version clearly reveals that, for certain instances, the solver converges to local optima (both with flat and data-boosted initializations), as evidenced by the high maximum \texttt{sub} values. The \texttt{SP-O} approach also reaches the time limit in both configurations and shows relatively large $\texttt{gap}$ values, particularly for the tuned version. In this latter case, a significant number of instances are infeasible, indicating that the solver fails to identify a feasible solution within the time limit, even when using the original bounds. Specifically, for \texttt{SP-O}, no feasible solution is found in 59 of the 100 instances, while this number is reduced to 50 and 39 instances for \texttt{SP-20} and \texttt{SP-100}, respectively. This shows that using tighter bounds derived from historical data helps the solver find feasible solutions more often within the time limit. For the standard 57-bus network, the data-boosted spatial-branching methods successfully reduce computational time with speedups of up to 654$\times$ over \texttt{SP-O}, but feasibility is fully preserved only when using $K=$ 100.

For the 118-bus networks, the interior-point methods continue to perform remarkably well, achieving low computational times and consistently reaching the global optimum for both initialization schemes. For this system size, however, some differences appear between formulations: the polar formulation is generally faster than the rectangular one, and the benefits of using a data-boosted initialization based on the nearest neighbor are particularly evident, as reflected in the reduced computational time when the interior-point method is applied to the rectangular formulation of the OPF problem with the data-driven initialization. The \texttt{SP-O} approach reaches the one-hour time limit but still attains the optimal solution, even though global optimality cannot be formally certified since $\texttt{gap}$ values are above the threshold. The data-boosted spatial-branching variants, on the other hand, are up to 321$\times$ faster while successfully delivering optimal solutions without any infeasible instances for $K=$ 100. This network thus represents the case where the data-boosted spatial-branching approaches exhibit their best performance. Nonetheless, even under these favorable conditions, their efficiency remains below that of the interior-point methods, regardless of whether local optima are present.

In summary, the results lead to several clear conclusions. For the standard configurations of all four networks, interior-point methods consistently find the global optimum in less than 10 seconds for all test instances, with the polar formulation combined with data-driven initialization yielding the best results for the largest network size. In contrast, even for smaller systems, the \texttt{SP-O} approach frequently reaches the time limit. While data-boosted variants of spatial branching reduce computational time and improve feasibility, their performance still falls short of that achieved by interior-point methods. In the tuned configurations designed to induce local optima, interior-point methods continue to find the global optimum for the 14-, 39-, and 118-bus networks, with local optima observed only in some instances of the 57-bus network. In these more challenging cases, data-boosted spatial branching also struggles, yielding a notable number of infeasible instances due to overly restricted bounds that exclude the true global optimum. Overall, across all network sizes and operating conditions, neither the standard nor data-boosted spatial-branching methods outperform interior-point approaches. This confirms that, despite promising advances, spatial branching for AC-OPF still demands further methodological and computational improvements before data-driven variants can match the efficiency of interior-point methods or the idealized performance of the benchmark \texttt{SP-$\epsilon$}.

\section{Conclusions} \label{sec:conclusions}

The AC Optimal Power Flow (AC-OPF) problem remains a central yet challenging task in power system operation due to its inherent non-convexity. While interior-point methods are widely used in practice for their computational efficiency, global optimization techniques such as spatial branching offer theoretical guarantees at a higher computational cost.

In this work, we introduced a data-boosted optimization framework that enhances both local and global solution strategies by leveraging historical operating data. Specifically, data-driven initialization improves the performance of interior-point methods, while learned variable bounds reduce the search space in spatial branching. This unified perspective enables a systematic evaluation of how learning can be used to accelerate fundamentally different optimization paradigms.

Our results show that modern interior-point solvers exhibit a high degree of practical robustness, consistently converging to globally optimal solutions across a wide range of test cases, including challenging instances with multiple local optima. Data-driven initialization further improves their computational performance, particularly for larger systems. In contrast, although data-boosted spatial branching significantly reduces computational effort, it remains substantially more demanding and may suffer from infeasibility when the learned bounds exclude the global optimum.

Overall, these findings highlight a clear gap between theoretical guarantees and practical performance in AC-OPF solution methods. While data-driven techniques can effectively enhance both approaches, they do not fundamentally alter this gap: interior-point methods remain the most efficient and reliable option in practice, whereas global optimization approaches still require substantial advances to become competitive at scale.

\printcredits

\section*{Declaration of competing interest}
The authors declare that they have no competing financial interests 
or personal relationships that could have appeared to influence the 
work reported in this paper. Ignacio Repiso, Salvador Pineda, and 
Juan Miguel Morales received financial support from the Spanish 
Ministry of Science and Innovation.

\section*{Acknowledgments}
This work was supported in part by the Spanish Ministry of Science and Innovation (AEI/10.13039/501100011033) through project PID2023-148291NB-I00 and by the Consejería de Universidad, Investigación e Innovación de la Junta de Andalucía through FEDER funds (grant TEP967-G-FEDER). I. Repiso, S. Pineda, and J. M. Morales are with the research group OASYS, University of Malaga, Malaga 29071, Spain. The work of I. Repiso was supported by the Spanish Ministry of Science and Innovation training program for PhDs with fellowship number by PRE2023-002227. Finally, the authors thankfully acknowledge the computer resources, technical expertise, and assistance provided by the SCBI (Supercomputing and Bioinformatics) center of the University of M\'alaga.

\section*{Data availability}
Data will be made available on request.

\bibliographystyle{cas-model2-names}

\bibliography{cas-refs}



\end{document}